\documentclass[11pt]{article}
\usepackage{amsmath,amsthm,amssymb}

\chardef\bslash=`\\ \hfuzz1pc
\newtheorem{thm}{Theorem}[section]
\newtheorem{prop}{Proposition}[section]
\newtheorem{rem}{Remark}[section]
\newtheorem{defn}{Definition}[section]
\numberwithin{equation}{section}

\begin{document}

\title{On extended umbral calculus, oscillator-like algebras and Generalize Clifford Algebra}
\author{A.K.Kwa\'sniewski
\\
Institute of Computer Science, Bia{\l}ystok University\\
PL-15-887 Bia{\l}ystok, ul.Sosnowa 64, POLAND\\
High School of Mathematics and Applied Computer Science,\\
PL-15-021 Bia{\l}ystok, ul. Kamienna 17, POLAND\\
}

\maketitle

\begin{abstract}
Some quantum algebras build from deformed oscillator algebras ala
Jordan-Schwinger may be described in terms of a particular case of
$psi$-calculus. We give here an example of a  specific relation
between such certain quantum algebras and generalized Clifford
algebras also in the context of  Levy -Leblond`s azimuthal
quantization of angular momentum which   was interpreted
afterwards as the finite dimensional quantum mechanics by
Santhanam et. all. $\psi$-calculus used for that as a framework is
that of classical operator calculus of Rota. By its nature
$\psi$-umbral calculus supplies a simple mathematical underpinning
for $\psi$-deformed quantum-like oscillator algebras and - at
least for the $\psi_{n}(q)=[n_{q}!]^{-1}$ case [1-3]. It provides
the natural underpinning for quantum group investigation.
Moreover- the other way around -one may formulate q-extended
finite operator calculus with help of the "quantum q-plane"
q-commuting variables $\ldots$ $\psi$-calculus is expected to be
useful in $C^{*}$ algebraic [4] description of "$\psi$-quantum
processes" with various parastatistics [5].

\end{abstract}

\section{\protect Extended Umbral Calculus in brief}
{\small The foundations of what we are going to call
"$\psi$-extension of Rota finite operator calculus" we owe  to
Viscov  [6,7]. $\psi$-extended umbral  calculus is arrived at
[6,7,8,9] by considering not only polynomial sequences of binomial
type but also of $\{s_{n}\}_{n\geqslant 1}$-binomial type where
$\{s_{n}\}_{n\geqslant 1}$-binomial coefficients are defined with
help of  the generalized factorial $n_{s}!=s_{1}s_{2}s_{3}\ldots
s_{n}$; $S=\{s_{n}\}_{n\geqslant 1}$ is an arbitrary $F$ valued
sequence with the condition $s_{n}\neq 0$, $n \in N$. F is a field
of char $F=0$ then the extensions rely on the notion of
$\partial_{\psi}$-shift invariance of $\partial_{\psi}$-delta
operators. Here $\partial_{\psi}$ denotes the $\psi$-derivative
i.e. $\partial_{\psi}x^{n}=n_{\psi}x^{n-1}$; $n\geqslant 0$; (then
$\partial_{\psi}$-linearly extended) and $n_{\psi}$ denotes the
$\psi$-deformed number where in conformity with Viscov notation
$n_{\psi}\equiv \psi_{n-1}(q)\psi_{n}^{-1}(q)$; $n_{\psi}!\equiv
\psi_{n}^{-1}(q)\equiv
n_{\psi}(n-1)_{\psi}(n-2)_{\psi}(n-3)_{\psi}\dots
2_{\psi}1_{\psi}$; $0_{\psi}!=1$. We choose to work with $\Im$-
the family of functions sequences such that $\Im=\{\psi;R\supset
[a,b];q\in[a,b];\psi(q):Z\rightarrow F;
\psi_{0}(q)=1;\psi_{n}(q)\neq 0;\psi_{-n}(q)=0; n\in N\}$. With
the choice $\psi_{n}(q)=[R(q^{n})!]^{-1}$ and
$R(x)=\frac{1-x}{1-q}$ we get the well known q-factorial
$n_{q}!=n_{q}(n-1)_{q}!$; $1_{q}!=0_{q}!=1$ while the
$\psi$-derivative $\partial_{\psi}$ becomes now the Jackson`s
derivative (see: [10-15]) $\partial_{q}$:
$(\partial_{q}\varphi)(x)=\frac{\varphi(x)-\varphi(qx)}{(1-q)x}$.
A polynomial sequence $\{p_{n}\}_{0}^{\infty}$ is then of
$\psi$-binomial type if it satisfies the recurrence
$$E^{y}(\partial_{\psi})p_{n}(x)\equiv p_{n}(x+_{\psi}y)=
\sum_{k\geqslant 0} {n\choose k}_{\psi}p_{k}(x)p_{n-k}(y)$$ where
${n \choose k}_\psi\equiv \frac
{n_{\psi}^{\underline{k}}}{k_{\psi}!}$. Here
$E^{y}(\partial_{\psi})\equiv
\exp_{\psi}\{y\partial_{\psi}\}=\sum_{k=0}^{\infty}\frac{y^{k}\partial_{\psi}^{k}}{n_{\psi}!}$
is a generalized translation operator and $\partial_{\psi}$-shift
invariance is defined accordingly. Namely  we work with
$\sum_{\psi}$ which is the algebra of $F$-linear operators acting
on the algebra $P$ of polynomials. These operators are
$\partial_{\psi}$-shift invariant operators $T$ i.e. $\forall
\alpha \in F; [T; E^{\alpha}(\partial_{\psi})]=0$. One then
introduces the notion of $\partial_{\psi}$-delta operator
according to
\begin{defn}
Let  $Q(\partial_{\psi}):P\rightarrow P$; the linear operator
$Q(\partial_{\psi})$ is a $\partial_{\psi}$-delta operator iff
 a)$Q(\partial_{\psi})$ is $\partial_{\psi}$-shift invariant;  b)$Q(\partial_{\psi})(id)=const\neq 0$.
 \end{defn}
  The strictly related notion is that of the $\partial_{\psi}$-basic polynomial sequence:
\begin{defn}
Let $Q(\partial_{\psi}):P\rightarrow P$; be the
$\partial_{\psi}$-delta operator. A polynomial sequence
$\{p_{n}\}_{n\geqslant 0}$; deg $p_{n}=n$ such that:
1)$p_{0}(x)=1$; 2)$p_{n}(0)=0; n\geqslant 1;$
3)$Q(\partial_{\psi})p_{n}=n_{\psi}p_{n-1}$ is called the
$\partial_{\psi}$-basic polynomial sequence of the
$\partial_{\psi}$-delta operator $Q(\partial \psi)$.
\end{defn}
After that and among many others the important Theorem (\ref{tw1})
might be proved using the fact that $\forall Q(\partial_{\psi})$
$\exists !$ invertible $S_{\partial_{\psi}}\in\sum_{\psi}$ such
that $Q(\partial_{\psi})=\partial_{\psi}S_{\partial_{\psi}}$.
\begin{thm}[$\psi$-Lagrange and $\psi$-Rodrigues formulas] \label{tw1}
Let $\{p_{n}(x)\}_{n=0}^{\infty}$ be $\partial_{\psi}$-basic
polynomial sequence of the $\partial_{\psi}$-delta operator
$Q(\partial_{\psi})$: Let $Q(\partial_{\psi})=\partial_{\psi}S$.
Then for $n>0$:
\begin{enumerate}
  \item $p_{n}(x)=Q(\partial_{\psi})'S^{-n-1}x^{n}$;
  \item
  $p_{n}(x)=S^{-n}x^{n}-\frac{n_{\psi}}{n}(S_{-n})'x^{n-1}$;
  \item $p_{n}(x)=\frac{n_{\psi}}{n}\hat{x}_{\psi}S^{-n}x^{n-1}$
  \item
  $p_{n}(x)=\frac{n_{\psi}}{n}{\hat{x}}_{\psi}(Q(\partial_{\psi})')^{-1}p_{n-1}(x)$ ($\psi$-Rodrigues formula).
\end{enumerate}

\end{thm}
Here we used the properties of the Pincherle $\psi$-derivative:

\begin{defn}(compare with  (17) in  [7])
The Pincherle $\psi$-derivative i.e. the linear map
$':\sum_{\psi}\rightarrow \sum_{\psi}$.
$$T'=T\hat{x}_{\psi}-\hat{x}_{\psi}T\equiv [T,\hat{x}_{\psi}]$$ where linear
operator $\hat{x}_{\psi}:P\rightarrow P$; is defined by
$$\hat{x}_{\psi}x^{n}=\frac{\psi_{n-1}(q)(n+1)}{\psi_{n}(q)}x^{n+1}=\frac{n+1}{(n+1)_{\psi}}x^{n+1};
n\geqslant 0.$$
\end{defn}

The  more general class constitute Sheffer
$\partial_{\psi}$-polynomials defined as:
\begin{defn}
A polynomial sequence $\{s_{n}(x)\}_{n=0}^{\infty}$ is called  the
sequence $\{s_{n}(x)\}_{n=0}^{\infty}$ of Sheffer
$\partial_{\psi}$-polynomials of the $\partial_{\psi}$-delta
operator iff (1)$s_{0}(x)=c\neq 0$, (2)
$Q(\partial_{\psi})s_{n}(x)=n_{\psi}s_{n-1}(x)$.
\end{defn}
The following proposition  relates Sheffer
$\partial_{\psi}$-polynomials of the $\partial_{\psi}$-delta
operator $Q(\partial_{\psi})$ to the unique
$\partial_{\psi}$-basic polynomial sequence of the
$\partial_{\psi}$-delta operator $Q(\partial_{\psi})$:
\begin{prop} \label{pr1}
Let $Q(\partial_{\psi})$ be a $\partial_{\psi}$-delta operator
with $\partial_{\psi}$-basic polynomial sequence
$\{q_{n}(x)\}_{n=0}^{\infty}$. Then $\{s_{n}(x)\}_{n=0}^{\infty}$
is a sequence of  Sheffer $q$-polynomials of the
$\partial_{\psi}$-delta operator $Q(\partial_{\psi})$ iff there
exists a $\partial_{\psi}$-shift invariant operator
$S_{\partial_{\psi}}$ such that $s_{n}(x)=S^{-1}q_n(x)$.
\end{prop}
 The family of Sheffer $\partial_{\psi}$-polynomials`
sequences $\{s_{n}(x)\}_{n=0}^{\infty}$ corresponding to the fixed
$\partial_{\psi}$-delta operator $Q(\partial_{\psi})$ is labeled
by elements of the abelian group of all $\partial_{\psi}$-shift
invariant invertible operators $S$. It is an orbit of this
group.\\
\textbf{Examples} According to Proposition \ref{pr1} with
$Q(\partial_{q})=\partial_{q}$ and
$S_{\partial_{q}}=E^{\alpha}(\partial_{q})\exp_{\psi}\{\partial_{q}^{2}\}$
we get q-Hermite polynomials while with choice
$Q(\partial_{q})=\frac{\partial_{q}}{\partial_{q}-1}$ and
$S_{\partial_{q}}=(1-\partial_{q})^{\alpha +1}$ we obtain
q-Laguerre polynomials $L_{n,q}^{(\alpha)}(x)$ of order $\alpha$.
So ($\psi$-extension of finite Rota calculus includes q-Hermite,
q-Laguerre polynomials $L_{n,q}^{(\alpha)}(x)$ of order $\alpha$
with their $\psi$-correspondents. These are already well known
q-Sheffer polynomials [16,9]. Specifically q-Laguerre polynomials
$L_{n,q}^{(-1)}(x)\equiv L_{n,q}(x)$ form the $\partial_{q}$-basic
polynomial sequence of $\{L_{n,q}(x)\}_{n\geqslant 0}$ the
$\partial_{q}$ operator
$$Q(\partial_{q})=-\sum_{k=0}^{\infty}\equiv \frac
{\partial_{q}}{\partial_{q}-1}\equiv
-[\partial_{q}+\partial_{q}^{2}+\partial_{q}^{3}+\dots].$$ Using
then Theorem (\ref{tw1}) one arrives at explicit form of
$L_{n,q}(x)$. Namely:\\
$L_{n,q}(x)=\frac{n_{q}}{n}\hat{x}_{q}(\frac{1}{\partial_{q}-1})^{-n}x^{n-1}
=\frac{n_{q}}{n}\hat{x}_{q}(\partial_{q}-1)^{n}x^{n-1}=
\frac{n_{q}}{n}\hat{x}_{q}\sum_{k=1}^{n}(-1)^{k}{n \choose
k}_{q}\partial_{q}^{n-k}x^{n-k}=\frac{n_{q}}{n}\sum_{k=1}^{n}(-1)^{k}{n
\choose
k}(n-k)_{q}^{\underline{n-k}}\frac{k}{k_{q}}x^{k}=\frac{n_{q}}{n}\sum_{k=1}^{n}(-1)^{k}\frac{n_{q}!}{k_{q}!}
\frac{(n-1)_{q}^{\underline{n-k}}}{(n-k)_{q}!}\frac{k}{k_{q}}x^{k}.$
\\So finally
\begin{equation}\label{e1}
L_{n,q}(x)=\frac{n_{q}}{n}\sum_{k=1}^{n}(-1)^k\frac
{n_{q}!}{k_{q}!}{{n-1}\choose {k-1}}_{q}\frac{k}{k_{q}}x^{k}
\end{equation}
Note: $\psi$-extended case is covered in this example just by
replacement $q \to \psi$. Let us also  stress here again  that
q-deformed quantum oscillator algebra provides a natural setting
for  q-Laguerre polynomials and q-Hermite polynomials [40,41].
$sl_{q}(2)$ and the q-oscillator algebra give rise to basic
geometric functions as matrix elements of certain operators in
analogy with Lie theory [19]  Also automorphisms of the
q-oscillator algebra lead to Sheffer q-polynomials  for example to
q-generalization of the Charlier polynomials [19].
\section{\protect q-extended quantum oscillator and Extended Umbral Calculus}
$\partial_{q}$-delta operators and their duals and similarly
$\partial_{\psi}$-delta operators with their duals   provide us
with pairs of generators of $\psi$-deformed quantum oscillator
algebras (see Remark 2.2)- possible candidates to describe
parastatistical behavior of some processes [5]. If
$\psi_{n}(q)=[R(q^{n})!]^{-1}$ and $R(x)=\frac{1-x}{1-q}$ then we
get the well known q-deformed oscillator dual pair of operators
leading to the corresponding $C^{*}$ algebra description of
q-Heisenberg-Weyl algebra. These oscillator-like algebras
generators and q-oscillator-like algebras generators are
encountered explicitly or implicitly in [1,2] and in many other
subsequent references - see [30,4,33] and references therein. In
many such references [32,33,19] q-Laguerre and q-Hermite or
q-Charlier polynomials appear which are just either Sheffer
$\psi$-polynomials or just $\partial_{\psi}$-basic polynomial
sequences of the $\partial_{\psi}$-delta operators
$Q(\partial_{\psi})$ for $\psi_{n}(q)=\frac{1}{R(q^{n})!}$;
$R(x)=\frac{1-x}{1-q}$ and corresponding choice of
$Q(\partial_{\psi})$ functions of $\partial_{\psi}$: $Q=id$. The
case $\psi_{n}(q)=\frac{1}{R(q^{n})!}$: $n_{\psi}=n_{k}$;
$\partial_{\psi}=\partial_{R}$ and $n_{\psi(q)}=n_{R(q)}=R(q^{n})$
appears implicitly in [21] where advanced theory of general
quantum coherent states is being developed. Among others also in
[23] it was noticed that commutation relations for the
q-oscillator-like algebras generators from [1 , 2] and others (see
[32 , 4]) might be chosen in appropriate operator variables to be
of the form [23]:

\begin{equation}\label{e21}
 AA^{+}-\mu A^{+}A=1; \mu=q^{2}
\end{equation}
As for the Fock space representation of normalized eigenstates
$|n\rangle$ of excitation number operator N various q-deformations
of the natural number n are used in literature on quantum groups
and at least some families of quantum groups  may be constructed
from q-analogues of Heisenberg  algebra  [1,2,23,3,34,35]. In
fact, these q-oscillator algebras generators are the so called
$\partial_{q}$-delta operators $Q(\partial_{q}$ i.e. basic objects
of the q-extended finite operator calculus of Rota. (Of course
$\partial_{q}\hat{x}-q\hat{x}\partial_{q}=id.$). The known
important fact is that the "q-Canonical Commutation Relations"
$AA^{+}-qA^{+}A=1$ lead [4] to the q-deformed spectrum of
excitation number operator N and to various parastatistics [5].
More possibilities result from considerations of  Wigner [36]
extended by the authors of [5]. We therefore hope that the
$\psi$-calculus of Rota to be developed here might be useful in a
$C^{*}$ algebraic [4] description of "$\psi$-quantum processes" -
if any - with various parastatistics [5]. Here in below via series
of definitions we shall propose a $\psi$-extension of the
q-oscillator model algebra using basic concepts of
$\psi$-extension of calculus of Rota.
\begin {defn} Let $\{p_{n}\}_{n\geqslant 0}$
be the $\partial_{q}$-basic polynomial sequence of the
$\partial_{q}$-delta operator $Q(\partial_{q})$. A linear map
$\hat{x}_{Q(\partial_{\psi})}: P \to P$;
$\hat{x}_{Q(\partial_{\psi})}p_{n}=p_{n+1}; n\geqslant 0$ is
called the operator dual to $Q(\partial_{q})$ .
\end{defn}
 For $Q=id$ we have: $\hat{x}_{Q(\partial_{q})}\equiv \hat{x}_{\partial_{q}}\equiv \hat{x}.$
 \\
\textbf{Comment:} Dual in the above sense corresponds to adjoint
in q-umbral calculus language of linear functionals' umbral
algebra (see : Proposition 1.1.21 in [16]).

\begin{defn}
Let $\{p_{n}\}_{n\geqslant 0}$ be the $\partial_{\psi}$-basic
polynomial sequence of the $\partial_{\psi}$-delta operator
$Q(\partial_{\psi})=Q$ then the $\hat{q}_{\psi Q}$-operator is a
liner map; $\hat{q}_{\psi,Q}:P \to P$;
$\hat{q}_{\psi,Q}p_{n}=\frac{(n+1)_{\psi}-1}{n_{\psi}}p_{n}$;
$n\eqslantgtr 0$.
\end{defn}
We call the $\hat{q}_{\psi,Q}$ operator the
$\hat{q}_{\psi,Q}$-mutator operator.\\
 \textbf{Note}: For $Q=id$, $Q(\partial_{\psi})=\partial_{\psi}$
the natural notation is $\hat{q}_{\psi,id}\equiv \hat{q}_{\psi}$.
For $Q=id$ and $\psi_{n}(q)=\frac{1}{R(q^{n})!}$ and
$R(x)+\frac{1-x}{1-q}\hat{q}_{\psi,Q}\equiv \hat{q}_{R,id}\equiv
\hat{q}_{R} \equiv \hat{q}_{q,id}\equiv \hat{q}_{q}\equiv \hat{q}$
and $\hat{q}_{\psi,Q}x^{n}=q^{n}x^{n}$.
\begin{defn}
Let $A$ and $B$ be linear operators acting on $P$; $A:P \to P$;
$B:P\to P$. Then $AB-\hat{q}_{\psi,Q}BA\equiv
[A,B]_{\hat{q}_{\psi,Q}}$ is called $\hat{q}_{\psi,Q}$-mutator of
$A$ and $B$ operators.
\end{defn}
\textbf{Note:}
$Q(\partial_{\psi})\hat{x}_{Q(\partial_{\psi})}-\hat{q}_{\psi,Q}\hat{x}_{Q(\partial_{\psi})}Q(\partial_{\psi})
\equiv
[Q(\partial_{\psi}),\hat{x}_{Q(\partial_{\psi})}]_{\hat{q}_{\psi,Q}}=id.$
This is easily verified in the $\partial_{\psi}$-basic
$\{p_{n}\}_{n\geqslant 0}$ of the $\partial_{\psi}$-delta operator
$Q(\partial_{\psi})$.
\\ Equipped with pair of operators $(Q(\partial_{\psi}),\hat{x}_{Q(\partial_{\psi})})$ and
$\hat{q}_{\psi,Q}$-mutator we have at our disposal all possible
representants of "$\psi$-canonical pairs" of differential
operators on the $P$ algebra. The meaning of the adjective:
"$\psi$-canonical" is explained by the content of the remark 2.2.
For important historical reasons however here is at first:
\begin{rem} The $\psi$-derivative is a particular example of a linear operator that
reduces by one the degree of any polynomial. In 1901 it was proved
[26] that every linear operator $T$ mapping $P$ into $P$ may be
represented as infinite series in operators $\hat{x}$ and $D$. In
1986 the authors of [27] supplied the explicit expression for such
series in most general case of polynomials in one variable; namely
according to the Proposition 1 from [27] one has: "Let $\triangle$
be a linear operator that reduces by one each polynomial. Let
$\{q_{n}(\hat{x})\}_{n\geqslant 0}$ be an arbitrary sequence of
polynomials in the operator $\hat{x}$. Then $T=\sum_{n\geqslant
0}q_{n}(\hat{x})\triangle ^{n}$ defines a linear operator that
maps polynomials into polynomials. Conversely, if $T$ is linear
operator that maps polynomials into polynomials then there exists
a unique expansion of the form
\begin{equation}\label{e22}
  T=\sum_{n\geqslant 0}q_{n}(\hat{x})\triangle ^{n}"
\end{equation}

\end{rem}
\textbf{Note:}In 1996 this was extended to algebra of many
variables polynomials [28].
\begin{rem}
The importance of the pair of dual operators: $Q(\partial_{\psi})$
and $\hat{x}_{Q(\partial_{\psi})}$ is reflected by the facts:
\\a) $Q(\partial_{\psi})\hat{x}_{Q(\partial_{\psi})}-\hat{q}_{\psi,Q}\hat{x}_{Q(\partial_{\psi})}Q(\partial_{\psi})
\equiv
[Q(\partial_{\psi}),\hat{x}_{Q(\partial_{\psi})}]_{\hat{q}_{\psi,Q}}=id.$
\\ b) Let $\{q_{n}(\hat{x}_{Q(\partial_{\psi})})\}_{n\geqslant 0}$ be an arbitrary sequence of polynomials in the
operator $\hat{x}_{Q(\partial_{\psi})}$. Then $T=\sum_{n\geqslant
0}q_{n}(\hat{x}_{Q(\partial_{\psi})}Q(\partial_{\psi})^{n}$
defines a linear operator that maps polynomials into polynomials.
Conversely, if $T$ is linear operator that maps polynomials into
polynomials then there exists a unique expansion of the form
\begin{equation}\label{e23}
  T=\sum_{n\geqslant
0}q_{n}(\hat{x}_{Q(\partial_{\psi})})Q(\partial_{\psi})^{n}
\end{equation}
\end{rem}
Equipped with pair of operators
$Q(\partial_{\psi}),\hat{x}_{Q(\partial_{\psi})}$ and
$\hat{q}_{\psi,Q}$-mutator we have at our disposal all possible
representants of "$\psi$-canonical pairs" of linear operators on
the $P$ algebra such that: a) the above  unique expansion
$T=\sum_{n\geqslant
0}q_{n}(\hat{x}_{Q(\partial_{\psi})})Q(\partial_{\psi})^{n}$ holds
b) we have the structure of $\psi$-umbral  or $\psi$-extended
finite operator calculus - coworking.
\section{\protect No $\psi$-analogue of quantum q-plane
formulation? [37,38]} In [16] Cigler and then  Kirchenhofer
defined the polynomial sequence $\{p_{n}\}_{0}^{\infty}$ of
q-binomial type equivalently by
\begin{equation}\label{e31}
  p_{n}(A+B)\equiv \sum_{k\geqslant 0}{n \choose
  k}_{q}p_{k}(A)p_{n-k}(B)
\end{equation}
where $[B,A]_{q}\equiv BA-qAB=0$. $A$ and $B$ might be interpreted
then as coordinates on quantum q-plane (see Ref.9 Chapter 4). For
example $A=\hat{x}$ and $B=y\hat{Q}$ where
$\hat{Q}\varphi(x)=\varphi(qx)$. If so then the following
identification holds: $$p_{n}(x+_{q}y)\equiv
E_{y}(\partial_{q})p_{n}(x)=\sum_{k\geqslant 0}{n\choose
k}_{q}p_{k}(x)p_{n-k}(y)=p_{n}(\hat{x}+y\hat{Q})\textbf{1}$$ Also
q-Sheffer polynomials $\{s_{n}(x)\}_{n=0}^{\infty}$ are defined
equivalently (see: 2.1.1. in [13]) by
\begin{equation}\label{e32}
  s_{n}(A+B)\equiv\sum _{k\geqslant 0}{n\choose
  k}_{q}s_{k}(A)p_{n-k}(B)
\end{equation}
where $[B,A]_{q}\equiv BA-qAB=0$ and $\{p_{n}(x)\}_{n=0}^{\infty}$
of q-binomial type. For example $A=\hat{x}$ and $B=y\hat{Q}$ where
$\hat{Q}\varphi(x)=\varphi(qx)$. Then the following identification
takes place:
\begin{equation}\label{e33}
  s_{n}(x+_{q}y)\equiv
  E^{y}(\partial_{q})s_{n}(x)=\sum_{k\geqslant0}{n \choose
  k}_{q}s_{k}(x)p_{n-k}(y)=s_{n}(\hat{x}+y\hat{Q})\textbf{1}
\end{equation}
This means that  one may formulate q-extended finite operator
calculus with help of the "quantum q-plane" q-commuting variables
$A,B:AB-qBA\equiv [A,B]_{q}=0$. One may now be tempted  - perhaps
in vain - to formulate  the basic notions of $\psi$-extended
finite operator calculus with help of the "quantum $\psi$-plane"
$\hat{q}_{\psi,Q}$-commuting variables
$A,B:[A,B]_{\hat{q}_{\psi,Q}}=0$ exactly in the same way as in
[16]. For that to try consider appropriate generalization of
$A=\hat{x}$ and $B=y\hat{Q}$ where this time the action of
$\hat{Q}$ on $\{x^{n}\}_{0}^{\infty}$ is to be found from the
condition $AB-\hat{q}_{\psi}BA\equiv [A,B]_{\hat{q}_{\psi}}=0$.
Acting with $[A,B]_{\hat{q}_{\psi}}$ on $\{x^{n}\}_{0}^{\infty}$
due to $\hat{q}_{\psi}x^{n}=\frac{(n+1)_{\psi}-1}{n_{\psi}}x^{n}$;
$n\geqslant 0$ one easily sees that now $\hat{Q}x^{n}=b_{n}x^{n}$
where $b_{0}=0$ and
$b_{n}=\prod_{k=1}^{n}\frac{(k+1)_{\psi}-1}{k_{\psi}}$ for $n>0$
is the solution of the difference equation:
$b_{n}-b{n-1}\frac{(n+1)_{\psi}-1}{n_{\psi}}=0$; $n>0$. With all
above taken into account one immediately verifies that for our $A$
and $B$ $\hat{q}_{\psi}$-commuting variables already
\begin{equation}\label{e34}
  (A+B)^{n}\neq \sum_{k\geqslant 0}{n \choose k}_{\psi}A^{k}B^{n-k}
\end{equation}
unless $\psi_{n}(q)=\frac{1}{R(q^{n})!}$; $R(x)=\frac{1-x}{1-q}$
hence $\hat{q}_{\psi,Q}\equiv \hat{q}_{R,id}\equiv
\hat{q}_{R}\equiv \hat{q}_{q,id}\equiv \hat{q}_{q}\equiv \hat{q}$
and $\hat{q}_{\psi,Q}x^{n}=q^{n}x^{n}$. Concluding: the above
identifications of polynomial sequence $\{p_{n}\}_{0}^{\infty}$ of
q-binomial type and Sheffer q-polynomials
$\{s_{n}(x)\}_{n=0}^{\infty}$ fail to be extended to the more
general $\psi$-case. This means that we can not formulate that way
the $\psi$-extended finite operator calculus with help of the
"quantum $\psi$-plane" $\hat{q}_{\psi,Q}$-commuting variables
$A,B:AB-\hat{q}_{\psi,Q}BA\equiv [A,B]_{\hat{q}_{\psi,Q}}=0$ while
considering algebra of polynomials $P$ over the field  $F$.
\section{\protect Relation to Quantum groups - Polar decomposition of $SUq(2;C)$ group algebra}
The standard basis of Lie algebra $su(2)=so(3)$ is well known to
be represented by:$$J_{3}=\sum_{m=-j}^{j}m|jm\rangle;$$
\begin{equation}\label{e41}
  J_{+}=\sum_{m=-j}^{j}\sqrt{(j-m)(j+m+1)}|j(m+1)\rangle\langle jm|;
\end{equation}
$$J_{-}=\sum_{m=-j}^{j}\sqrt{(j+m)(j-m+1)}|j(m+1)\rangle\langle jm|$$
In [1] Biedenharn  proposed a new realization of quantum group
$SU_{q}(2)$ and in order to realize generators of a q-deformation
$U_{q}(su(2))$ of the universal enveloping algebra of the Lie
algebra $su(2)$ he defined a pair of mutual commuting q-harmonic
oscillator systems (al\`{a} Jordan-Schwinger approach to $su(2)$
Lie algebra). At the same time [2] Mac Farlane had also discovered
q-oscillator description of $SU_{q}(2)$ (al\`{a} Jordan-Schwinger
approach to $su(2)$ Lie algebra). The generators of a
q-deformation $U_{q}(su(2))$ of the universal enveloping algebra
of the  Lie algebra $su(2)$ (called by physicists "generators of
the quantum group $SU_{q}(2)$" - which is not a group!) are
obtained from (4.1) by one of several possible q-deformations. In
Biedenharn's and Mac Farlane's case one uses the following
q-deformation of numbers and operators:
\begin{equation}\label{e42}
  [x]_{q}=\frac{q^{x}-q^{-x}}{q-q^{-1}}.
\end{equation}
Thus  q-deformed \ref{e41} representation of  generators now reads
 $$J_{3}=\sum_{m=-j}^{j}m|j,m\rangle_{q};$$
\begin{equation}\label{e43}
 J_{+}=\sum_{m=-j}^{j}\sqrt{[j-m]_{q}[j+m+1]_{q}}|j,(m+1)\rangle_{q  q}\langle
 j,m|;
\end{equation}
$$J_{-}=\sum_{m=-j}^{j}\sqrt{[j+m]_{q}[j-m+1]_{q}}|j,(m+1)\rangle _{q  q}\langle
j,m|,$$ where
\begin{equation}\label{e44}
  |j,m\rangle
  _{q}=|j+m\rangle_{q}|j-m\rangle_{q}=\frac{{a^{+}_{1q}}^{j+m}{a^{+}_{2q}}^{j-m}}{[j+m]_{q}![j-m]_{q}!}|0\rangle_{q}
\end{equation}
and $a^+_{1q}$, $a^+_{2q}$ represent two mutually commuting
creation operators of  q-quantum  oscillators. Corresponding
commutation relations of the generators of a q-deformation
$U_{q}(su(2))$ of the universal enveloping algebra of the Lie
algebra  $su(2)$ are of the familiar though now q-deformed form
[1],[2],[29]:
\begin{equation}\label{e45}
  [J_{3},J_{+}]=J_{+};  [J_{3},J_{-}]=-J_{-};
  [J_{+},J_{-}]=[2J_{3}]_{q}.
\end{equation}
From (\ref{e43}) one may derive [29] the polar decomposition of
the generators $J_{+}$, $J_{-}$:
\begin{equation}\label{e46}
 J_{+}=\sqrt{J_{+}J_-}\sigma_1^{-1}=\sigma_1^{-1}\sqrt{J_-J_+}\qquad
 J_-=\sqrt{J_+J_-}\sigma_1=\sigma_1\sqrt{J_-J_+},
\end{equation}
where $\sigma_1$ is the first of the two generators of generalized
Pauli algebra [30]
\begin{equation}\label{e47}
  \sigma_1=\left(\begin{array}{cccccccc}
    0 & 1 & 0 & 0 & \dots & 0 & 0 \\
    0 & 0 & 1 & 0 & \dots & 0 & 0 \\
    0 & 0 & 0 & 1 &  \dots & 0 & 0 \\
    . & . & . & . & . & . & . \\
    0 & 0 & 0 & 0 & \dots & 1 & 0 \\
    0 & 0 & 0 & 0 & \dots & 0 & 1 \\
    1 & 0 & 0 & 0 & \dots & 0 & 0 \
  \end{array}\right)=(n\times n)
\end{equation}
The second generator $\sigma_2$ has been also used in [29] in
order to remark on relevance of such a pair $\sigma_1$, $\sigma_2$
to $GL_\omega(2;C)$ properties
$(q\equiv\omega\equiv\exp\{\frac{2\pi i}{n}\})$. It is to be noted
here that the polar decomposition for undeformed $su(2)=so(3)$
algebra of undeformed quantum angular momentum had been performed
already by Levy -Leblond in [31]. There he had interpreted such a
polar decomposition as the "azimuthal quantization of angular
momentum". Following the authors of  [29] it should be noted here
that \textit{generalized Pauli algebra  appears in q-deformed and
in undeformed case of polar decomposition in the same way (4.6)}.
Neither  Levy -Leblond nor the authors of  [29] had realized that
they are dealing with generalized Pauli algebra [30]. These has
been realized afterwards by T.S. Santhanam  [32]. \\ In the
notation of [30] and papers quoted there
\begin{equation}\label{e48}
  \sigma_2\equiv U=\omega^Q=\exp\{\frac{2\pi i}{n}Q\}=
  \left(\begin{array}{cccc}
    1 & 0 & \dots & 0 \\
    0 & \omega & \dots & 0 \\
 . & . & \dots & . \\
    0 & 0 & \dots & \omega^{n-1} \
  \end{array}\right)
\end{equation}
where
\begin{equation}\label{e49}
  Q=\left(\begin{array}{cccc}
    1 & 0 & \dots & 0 \\
    0 & 1 & \dots & 0 \\
 . & . & \dots & . \\
    0 & 0 & \dots & n-1 \
  \end{array}\right)
\end{equation}
\begin{equation}\label{e410}
  \sigma_1\equiv V=\omega^P=\exp\{\frac{2\pi i}{n}P\}=
  \left(\begin{array}{ccccc}
    0 & 1 & 0 & \dots & 0 \\
    0 & 0 & 1 & \dots & 0 \\
    0 & 0 & 0 & \dots & 0 \\
 . & . & .& \dots & 1 \\
    1 & 0 & 0 & \dots & 0 \
  \end{array}\right)
\end{equation}
where
\begin{equation}\label{e411}
  P=S^\dagger QS=(P_{\alpha,\kappa}),
\end{equation}
$$P_{\alpha,\kappa}= \left\{
\begin{array}{ll}
  0 & \alpha=\kappa \\
  \lbrack \overline{\omega}^{\alpha-\kappa}-1\rbrack^{-1} & \alpha\neq\kappa
\end{array}
\right.$$ and
\begin{equation}\label{e412}
  s=(\langle \widetilde{k}|l\rangle)=\frac{1}{\sqrt{n}}(\omega^{kl})_{k,l\in Z_n}
\end{equation}
is the Sylvester matrix. Formulas (\ref{e48})-(\ref{e412}) contain
the main information on quantum kinematics of the finite
dimensional quantum mechanics as here we interpret polar
decomposition of quantum angular momentum algebra $su(2)=so(3)$
formalism as a model of  finite dimensional quantum mechanics with
the classical phase space being the torus $\textbf{Z}_n\times
\textbf{Z}_n$(see [30] and references therein). This possibility
was already considered by Weyl in [33]. "Azimuthal quantization of
angular momentum" was interpreted afterwards as the finite
dimensional quantum mechanics by Santhanam et. all [32]. The
considerations of this section allow us to hope to elaborate soon
more on the q-deformed finite dimensional quantum mechanics
treated as an interpretation of q-deformed  $su(2)$ algebra of
q-deformed angular momentum [34,35]. A specific relation between
certain quantum algebras and generalized Clifford algebras was
recently discovered in [36]. More elaborated presentation of the
content of  Sections I, II and III  one may find in [37,38].

\end{document}